\documentclass[11pt]{amsart}
\usepackage{amsbsy,amssymb,amsthm,amsmath}
\usepackage{syntonly}

\newcommand\what{\widehat}

\newcommand\hY{{\what Y}}

\newcommand\R{\mathbb R}

\newcommand\order{\text{order}}
\renewcommand\Re{\text{Re}}

\newcommand\half{{1 \over 2}}

\newcommand\Tr{\rm tr}

\newcommand\ep{\varepsilon}

\newcommand\vf{{\bf f}}

\newcommand{\re}[1]{(\ref{#1})}

\newcommand\normalpar{{\parindent .5cm}}
\renewcommand\proof{\noindent {\bf Proof.} \quad}
\renewcommand\endproof{{\hfill $\square$}}

\theoremstyle{plain}
\newtheorem{thm}{Theorem}[section]
\newtheorem{lem}[thm]{Lemma}
\newtheorem{defi}[thm]{Definition}

\newtheorem{cor}[thm]{Corollary}
\newtheorem{pro}[thm]{Proposition}

\newtheorem{rmk}[thm]{Remark}

\catcode`\@=11
\@addtoreset{equation}{section}   

\normalpar
\newcommand\Purdue{Department of
Mathematics, Purdue University,\\ West Lafayette, IN 47907, USA.}
\newcommand\Amiens{LAMFA CNRS UMR 6140 \\
Universit\'e de Picardie Jules Verne, \\
 33, rue Saint-Leu, 80039 Amiens, France.}
\title[09/16/04]{Long-Time Asymptotic Behavior of Dissipative
Boussinesq System}

\author[M.\ Chen \& O. Goubet] {M. Chen$^{1}$ and  O. Goubet$^{2}$}
\dedicatory{\vspace{-10pt}\normalsize{$^1$\Purdue \\
$^{2}$\Amiens}}

\begin{document}

\begin{abstract} In this paper, we study various dissipative mechanics
associated with the Boussinesq systems which model  two-dimensional small
amplitude
long wavelength water waves.
We will show that the decay rate for the damped one-directional model
equations, such as the KdV and BBM equations, holds for some of the damped
Boussinesq systems which model two-directional waves. 

\end{abstract}


\maketitle

\section{introduction}

Considered here are  waves on the surface of an inviscid fluid in a flat channel.
When one is  interested in the propagation of one-directional
irrotational 
small amplitude long waves, it is classical to model the waves  by the well-known
KdV (Korteweg-de Vries) equation (see \cite{Whitham74})
\begin{equation*}\label{kdv}
u_t+u_x+u_{xxx}+uu_x=0,
\end{equation*}
or its regularized version, the so-called regularized
long wave equation or BBM (Benjamin-Bona-Mahony) equation,
\begin{equation*}\label{bbm}
u_t+u_x-u_{txx}+uu_x=0.
\end{equation*}
When one is dealing with two-directional waves, and  the effects of wave  interactions and/or
wave reflections  are not excluded from the study, a restricted
four-parameter family of systems (see \cite{BCS02}), 
\begin{equation}
\label{abcd}
\begin{split}
&\eta_t+u_x+(u\eta)_x+au_{xxx}-b\eta_{xxt}=0, \\
&u_t+\eta_x+uu_x+c\eta_{xxx}-du_{xxt}=0,
\end{split}
\end{equation}
may be used. 
The
dimensionless variables $\eta(x,t)$, $u(x,t)$, $x$,  and $t$ are
scaled by the length scale $h_0$ and time scale $(h_0/g)^\half$
where $h_0$ denotes
the still water depth and $g$ denotes the acceleration of gravity. 
The variable $\eta(x,t)$ is the non-dimensional
deviation of the water surface from its undisturbed position and
$u(x,t)$ is the non-dimensional horizontal velocity at a height
above the bottom of the channel corresponding to $\theta h_0$ with
$0 \leq \theta \leq 1$. The constants $a, b, c, d$ are called dispersive constants
which satisfy the physical relevant
constraints
\begin{equation*}
(C0) \qquad a+b+c+d=\frac13 \quad \text{and} \quad c+d=\frac12 (1-\theta^2) \geq 0.
\end{equation*}

This class of systems contains some of the well-known systems, such as
the classical Boussinesq system
$
(a=b=c=0,\ d=1/3)
$
(see for example \cite{Boussinesq1871a,Peregrine72,TenWu94,Amick84,Schonbek81})
and
the Bona-Smith system $(a=0,\ b=d>0,\ c<0)$ \cite{BonSmi76}. It is
shown in \cite{BCS04} that
a physically relevant system in  \re{abcd} is linearly well posed in certain natural 
  Sobolev spaces if
the constants $a, b, c, d$ satisfy 
$$
(C1)\qquad b\geq 0,\ d\geq 0,\ a\leq 0,\ c\leq 0,
$$
or
$$
(C2)\qquad b\geq 0,\ d\geq 0,\ a=c>0.
$$
It is also shown in \cite{BonChe98,BPS81} that above systems have the capacity of 
capture  the main characteristics  of the flow in an idea fluid. But when the damping
effect is  comparable with the effects of nonlinearity and/or
dispersion,  as occurs in the  real laboratory-scale experiments
and in the fields (see \cite{BPS81,MahPri80,MeiLiu73,Miles67}),  it should  be considered  
in order for the model and its numerical results to  correspond in
detail with the experiments.
The full system
would be the Navier-Stokes equations with a free boundary,
which is  very difficult to handle both  theoretically and
numerically
(cf. \cite{TPL90,BHL96}).
Therefore, it is useful to construct
simpler model systems which are capable of  capturing  the main
properties of water waves under various special circumstances.

For example, many researchers have studied the
dissipative one-way propagation
model equations, such as the dissipative KdV and dissipative
regularized long-wave equations and their generalizations. As a
model to our study, we recall the results from \cite{ABS89} for
the dissipative BBM equation,
\begin{equation*}\begin{split}
u_t-u_{xxt}-\nu u_{xx}+uu_x=0, \\
u(x,0)=u_0(x)
        \end{split}\label{1.2}
\end{equation*}
where $\nu$ is a positive constant.
\begin{thm}\label{1} Assume $u_0$ is
in $L^1(\R)\cap L^2(\R)$, then there exists a constant $C$ such
that
\begin{equation}
\|u(t)\|_{L^2}\leq C (1+t)^{-1/4}. \label{1.3}
\end{equation}
\end{thm}

\noindent Here $L^1(\R)$ and $L^2(\R)$ are the classical Banach
spaces. A similar result holds for the corresponding dissipative KdV equation.
See also \cite{BPS81}, \cite{Karch00} and the references
therein. 

In this article, we aim to analyze the effect of dissipation on
systems \re{abcd} and study the decay rates of solutions $(\eta,u)$
toward zero.
We will restrict our  study to the cases where constants $a, b, c, d$
satisfy (C0)-(C1) or (C0)-(C2).
The goal of   this research  is to find the appropriate
dissipative
term (or terms) which will provide the
right amount of  energy dissipation for all wave numbers
while
keeping  the mass conserved. 

In this article, two kinds of
dissipations will be considered:

\noindent {\it Complete dissipation}: replacing the
$\begin{pmatrix} 0\\ 0\end{pmatrix}$ in the right-hand side of
\re{abcd} by  the vector $\begin{pmatrix} \eta_{xx}\\
u_{xx}\end{pmatrix}$, and 

\noindent
{\it Partial dissipation}: replace the $\begin{pmatrix} 0\\
  0\end{pmatrix}$
in the right-hand side of \re{abcd} by the vector $\begin{pmatrix}
0\\ u_{xx}\end{pmatrix}$.

We shall first study the decay rates of  solutions to
the linearized systems supplemented with either
complete or partial dissipations. These equations read
\begin{equation}
\begin{split}
&\eta_t+u_x+au_{xxx}-b\eta_{xxt}=\nu\eta_{xx},\\
&u_t+\eta_x+c\eta_{xxx}-du_{xxt}=u_{xx}, \end{split}\label{1.6}
\end{equation}
with $\nu=0$ or $\nu=1$. Systems
 which satisfy the {\it dichotomy} property in the Fourier space:
\begin{itemize}
\item decay as $t^{-1/4}$ for low frequencies (small $\xi$); 
\item decay as
$\exp(-\beta t)$ or $\exp(-\beta \xi^2 t)$ for high  frequencies (large $\xi$); 
\end{itemize}
will be identified and  studied.
It is shown in Section \ref{sec3} that the  dichotomy property 
 will lead to the decay rate $t^{-\frac14}$ for
$\|(\eta,u)\|_{L^2\times H^h}$
 ($h$ will be specified later). 

For later use, we shall emphasize (and it is easy to check) 
 that the dichotomy property
holds true for two fundamentally different equations: the 
linearized  BBM-Burgers equation $u_t-u_{xxt}+u_x-u_{xx}=0$ and the 
linearized KdV-Burgers equation
$u_t+u_{xxx}+u_x-u_{xx}=0$.  If
the high frequency part of the solution to a system  is  damped as exp$(-\beta
t)$, we say that the system belongs to the BBM-Burgers class since
solutions to linearized  BBM-Burgers equation feature this property. If the high
frequency part of the solution is damped as exp$(-\beta \xi^2t)$,
which is the case for linearized KdV-Burgers equation, 
then we say that the system belongs to
the KdV-Burgers class. The low frequency parts of the solutions to 
the linearized BBM-Burgers and KdV-Burgers equations behave
in a similar fashion.

The  main result in Section \ref{sec3}  is to classify the linearized systems according to this
property and to prove that for systems which satisfy the dichotomy property,
a decay rate comparing to \re{1.3} is valid. We shall also
present  some systems where  the decay rates can be
arbitrarily small, behaving as the solution of
$$
u_t-u_{xxt}+u=0, \quad u(x,0)=u_0(x), 
$$
where by Fourier transform,
$$ \what  u(t, \xi)=e^{-\frac{t}{1+\xi^2}} \what  u_0(\xi),$$
and therefore $\| u(.,t) \|_{L^2}$ could  decay arbitrarily slow.

In Section \ref{sec4}, we extend the linear theory  
to  nonlinear systems and show that
the decay rate
as \re{1.3} is valid for weakly dispersive systems, i.e.
systems with  $b>0$ and $d>0$,  and for some systems in the  KdV-Burgers
class with total dissipation, which include the 
KdV-KdV system ($a=c=\frac16, b=d=0$), with small initial data.
In Section \ref{sec5}, the decay rate with respect
 to $L^\infty$-norm is presented and in
Section \ref{sec6}, spectral method is used on several systems to demonstrate that 
the rates obtained in Section \ref{sec4} and Section
\ref{sec5}
 are sharp and the
constants involved in the bounds are reasonably sized. 

It is worth to note that there are other methods, 
such as the energy methods (like the
so-called Schonbek's splitting method
 applied
to
the classical Boussinesq system \cite{RSW99} in large dimensions),
can be used
in proving decay rate for solutions of these systems. This line of
study will be carried elsewhere. 
We believe that those methods will be helpful especially in the cases when $b=d$ so
some Hamiltonian is conserved (see \cite{BCS04}).

On the other hand, to remove the smallness assumption on the
initial data, the authors  in \cite{ABS89} used a kind of
Cole--Hopf transformation (that is valid for Burgers equation)
 and were able to control the extra-terms. We do not know
if this {\it tour de force} (feat of skill) is possible for the
systems in \re{abcd} with dissipation.

We complete this introduction by introducing some notations. Throughout the paper, the
standard notation on Sobolev spaces will be used. The $L^p(\R)$
norm  will be denoted as  $\|\cdot \|_{L^p}$ for $1 \leq p \leq
\infty$ and the $H^s$ norm will be denoted as  $\|\cdot \|_{H^s}$.
When several variables are involved, we may also set $L^p_x$ for
$L^p(\R)$ to specify that we compute the norm with respect to the
$x$-variable. The product space $X \times X$ will be abbreviated
by $X$ and a function $\vf=(f_1, f_2)$ in $X$ carries the norm
\begin{equation*}
\| \vf \|_{X}=\left(\| f_1 \|_X^2+\| f_2 \|_X^2\right)^{\half}.
\end{equation*}
The Euclidean norm of a vector is denoted by $|\cdot|$. We will use $C$ and $\beta$  as 
generic positive  constants whose values may change with each appearance.
Fourier transform of a function $f$ is denoted by either $\what f$
or $\mathcal {F}(f)$.

\section{Notations and Preparations}\label{sec2}

\subsection{ Some notations}

Consider $\nu\in\{0,1\}$.
As stated before, we plan to first estimate the decay rates of solutions to
the linear systems
\begin{equation}
\begin{split}
&\eta_t+u_x+au_{xxx}-b\eta_{xxt}=\nu\eta_{xx},\\
&u_t+\eta_x+c\eta_{xxx}-du_{xxt}=u_{xx},
\end{split}
\label{2.1}
\end{equation}
when $t$ goes to $+\infty$.

Following \cite{BCS04}, we introduce the Fourier multipliers
$$
\omega_1={1-a\xi^2\over 1+b\xi^2}\quad \text{ and
}\quad \omega_2={1-c\xi^2\over 1+d\xi^2}.
$$
Since $a, b, c, d$ satisfy  (C1) or (C2), $\omega_1\omega_2$ is
nonnegative and we denote
$$
\what H=\left({\omega_1\over \omega_2}\right)^{1/2}\text{ and }
\quad \sigma=(\omega_1 \omega_2)^{1/2},
$$
with the conventional notation
 ${0\over 0}=1$. We also denote

\def\ds{\displaystyle}

$$\alpha={\xi^2\over 1+b\xi^2} \quad \text{and} \quad \varepsilon=\ds{\xi^2\over 1+d\xi^2}.$$

\begin{rmk}
When a system satisfying  (C2) assumption is the subject of the study, 
$\omega_1$ and $\omega_2$ do change signs, but $\omega_1 \omega_2 \geq 0$.
\end{rmk}

\begin{defi}
Consider a nonnegative function $\xi \to \hat  \kappa(\xi)$. 
 The {\it order} of $\what  \kappa$ (when it exists) is defined as
the number $m$ such that
$$
\what  \kappa(\xi)\sim C|\xi|^m
$$
when $|\xi|\to +\infty$.
The (pseudo--differential) operator $\kappa$ with order $m$  is defined by setting
$$
\kappa u=v\quad \text{ iff } \quad \what  \kappa\what  u=\what  v.
$$
Therefore  $\kappa$ maps $L^2_x$ into
$H_x^{-\text{order}(\kappa)}$ (or $H^n_x$ into
$H_x^{n-\text{order}(\kappa)}$).
\end{defi}

Since \re{2.1} is a linear system, it is convenient to use the Fourier
transform.
Let $(\what \eta,\what  u)$ denote the Fourier transform of $(\eta, u)$ and set
$\hY= (\what \eta,\what  w)$ with  $\what  w=\what H\what  u$, then  \re{2.1} reads
\begin{equation}
\hY_t+A\hY=0,\label{2.2}
\end{equation}
where
$$A(\xi)=\begin{pmatrix} \nu \alpha&i\,{\rm sgn}(\omega_1)\xi\sigma\\
i\,{\rm sgn}(\omega_2)\xi\sigma&\varepsilon\end{pmatrix}$$ is the
{\it symbol}   of the linear (unbounded) operator  in 
\re{2.1}. Since we are dealing with a system,
$A(\xi)$ is a matrix.

By  multiplying  $\hY^*$  on \re{2.2}  and taking
the real part,
\begin{equation}
\begin{split}
& {1\over 2} {d\over dt}\int (|\what \eta(t,\xi)|^2+|\what w(t,\xi)|^2)d\xi\\
&+\nu\int\alpha(\xi)|\what \eta(t,\xi)|^2 d\xi
+\int\varepsilon(\xi)|\what w(t,\xi)|^2d\xi=0.
\end{split}
\label{2.3}
\end{equation}
Since $\alpha(\xi)$ and $\ep(\xi)$ are positive, 
\begin{equation}
E(t):=\int_{\R} |\hY(t,\xi)|^2 d\xi\label{2.4}
\end{equation}
decays towards $0$ as $t \rightarrow \infty$, where 
$|\hY|=(|\what \eta|^2+|\what w|^2)^{1/2}$
is the Euclidean norm on $\Bbb C^2$.

\subsection{Linear algebra}

We recall some facts from linear algebra and then
apply them  to the  dissipative systems \re{2.1}.

\begin{defi} Let  $M$ be a $2 \times 2$ matrix in the complex space, the norm of $M$ is defined by
$$
\|M\|=\sup_{Y\in\Bbb C^2\backslash\{0\}}\ {|M Y|\over |Y| }.
$$
\end{defi}

\begin{lem}\label{l2.3} Let  $\rho{(M)}$ denote the spectral radius of a matrix $M$
and $\Tr(M)$ denote the trace of $M$, then
$$
\|M\|=\rho{(M^* M)^{1/2}}\leq \Tr{(M^* M)}^{1/2}.
$$
\end{lem}

We are now going to bound  $E(t)$ (see \re{2.2} and \re{2.4}) by using the pointwise estimate
\begin{equation}
|\hY(t,\xi)|\leq \|e^{-t A}\| \ |\hY_0(\xi)|.\label{2.5}
\end{equation}

Noticing that the matrix $A$ can be written as $A=D+U$,  where
$D=\begin{pmatrix} \nu
\alpha&0\\ 0&\varepsilon\end{pmatrix}$ represents the  dissipation
terms and  $U=\begin{pmatrix} 0&i{\rm sgn}(\omega_1)\xi\sigma\\
i{\rm sgn}(\omega_2)\xi\sigma&0\end{pmatrix}$ is skew-symmetric.
When $D$ and $U$ commute,  the behavior of $\| e^{-tA}\|$ with respect
to $\xi$ is characterized
by the behaviors of $\ep$ and $\alpha$ via 
\begin{equation} \| e^{-tA}  \| \leq e^{-t \min \{ \nu\alpha(\xi), {\ep}(\xi) \}}.
 \label{rough2}
\end{equation}
But
when $D$ and $U$ do not commute,  more accurate estimate than
\re{rough2} can be obtained by studying
$e^{-tA}$ in detail.

We now  recall the following lemma (Theorem 9.28 from \cite{Finkbeiner66}).

\begin{lem}\label{lem2.4} There exists a unitary matrix $Q$
(i.e. $QQ^*=Q^*Q=I$)
 such that
$$
A=Q^*\begin{pmatrix} \lambda_1&z\\ 0&\lambda_2\end{pmatrix} Q,
$$
where $\lambda_1$ and $\lambda_2$ are the eigenvalues of $A$,
ordered by  $\Re(\lambda_1)\leq \Re(\lambda_2)$.
\end{lem}

As a consequence, one can prove (which is given at the end of this subsection)
\begin{lem} \label{lem2.5} There exists $C>0$ such that
\begin{equation}
\|\exp(-tA)\|\leq C\left(1+|z|\min \left(t, \frac{1}{|\lambda_2-\lambda_1|}\right)\right)\exp(-t\Re(\lambda_1)).\label{2.10}
\end{equation}
\end{lem}

It is easy to see that $\lambda_1$ and $\lambda_2$ are the roots of
the characteristic
equation
\begin{equation}
\lambda^2-\Tr(A) \lambda+\det(A)=0\label{2.7b}
\end{equation}
where
\begin{equation}
\Tr(A)=\lambda_1+\lambda_2=\nu \alpha +\ep \geq 0 \label{trA}
\end{equation}
and
\begin{equation}
\det(A)=\lambda_1\lambda_2=\nu \alpha \ep +\xi^2 \sigma^2 \geq 0. \label{detA}
\end{equation}

We now estimate
$\|\exp(-tA)\|$ by separating  the cases $\Delta\leq 0$ and
$\Delta > 0$ where $\Delta$ is the determinant of \re{2.7b}, namely
\begin{equation}
\Delta=\Tr(A)^2-4\det(A)=(\varepsilon-\nu\alpha)^2-4\xi^2\sigma^2.\label{2.7}
\end{equation}

\begin{lem}\label{lem2.6} 
For any $t>0$ and for any $\xi \in \R$, 
\item{$\bullet$} when  $\Delta\leq 0$ (perturbation range),
\begin{equation}
\|\exp (-tA)\|\leq C(1+\Tr(A)t)\exp\left(-{\Tr(A)\over 2}\ t\right)\leq C
\exp\left(-{\Tr(A)\over 4}\ t\right) ;\label{2.15}
\end{equation}
\item{$\bullet$} when $\Delta > 0$ (non-perturbation range),
\begin{equation}
\|\exp (-tA)\|\leq C\left(1+2|\xi| \sigma\min\left(t,
\frac{1}{\sqrt{\Delta}}\right)\right)\exp(-t \lambda_1),\label{2.16}
\end{equation}
where $\lambda_1$ satisfies
\begin{equation}
{\det(A)\over \Tr(A)}\leq\lambda_1\leq\min\left(\Tr(A),
{2\det(A)\over \Tr(A)}\right).\label{2.17}
\end{equation}
\end{lem}
\proof
It is worth to note from Lemma \ref{lem2.4} that
\begin{equation}
\Tr(A^* A)=|\lambda_1|^2 + |\lambda_2|^2 + |z|^2
=\nu^2\alpha^2+\varepsilon^2+2\xi^2\sigma^2.  \label{2.6}
\end{equation}

{\it When  $\Delta\leq 0$ (perturbation range):}
matrix $A$ has two conjugate complex eigenvalues $\lambda_1$ and $\lambda_2$ with
$$ |\lambda_1|=|\lambda_2|, \quad \Re(\lambda_1)=\Re(\lambda_2)=\frac{\Tr(A)}2, \quad |\lambda_1|^2=|\lambda_2|^2=\det(A).$$
Using \re{2.6} and then  \re{detA}-\re{trA}  leads to
\begin{equation*}\begin{split}
|z|^2&=\nu^2\alpha^2+\ep^2+2 \xi^2 \sigma^2 -2 |\lambda_1|^2\\
&=
\nu^2\alpha^2+\ep^2+2 \xi^2 \sigma^2 -2\det(A)=(\nu\alpha-\varepsilon)^2 \leq \Tr(A)^2.
        \end{split}
\label{2.8}
\end{equation*}
Hence \re{2.15} follows from  Lemma \ref{lem2.5}.

{\it When $\Delta > 0$ (non-perturbation range)}:
  $\Tr(A) \geq 0$  and $\det(A) \geq 0$ imply that  the matrix
 $A$ features two real eigenvalues $0\leq\lambda_1 < \lambda_2$.
 Then \re{2.6} leads to
\begin{equation*}
|z|^2=\nu^2\alpha^2+\ep^2+2 \xi^2 \sigma^2-\Tr(A)^2+2\det(A)=
4\xi^2\sigma^2\label{2.9}
\end{equation*}
and \re{2.16} is  proved by using
Lemma \ref{lem2.5}.
Since
$\lambda_2\leq \Tr(A)\leq 2\lambda_2$,
 one sees immediately that
\begin{equation*}
{\det(A)\over \Tr(A)}\leq\lambda_1=\frac{\det(A)}{\lambda_2}\leq
{2\det(A)\over \Tr(A)}
\end{equation*}
and  \re{2.17} follows.
\endproof

It is noted that when $\Delta\leq 0$, the dissipation can be
considered as a perturbation term with respect to the skew
symmetric operator. More precisely, the decay is the same as
pretending  $U$ and $D$ commute, up to a linear correction.

 When
$\Delta > 0$, this is no longer valid. In the first case, matrix $A$ 
has conjugate complex eigenvalues.
 In the latter case, $A$ has real positive eigenvalues
and the smallest one monitors the decay estimate.
  $\Delta=0$ is the
bifurcation point.

For the sake of
completeness, we now give the proof of Lemma \ref{lem2.5}.

\noindent
{\bf Proof of Lemma \ref{lem2.5}.} \quad
Straightforward computations lead to
\begin{equation*}
e^{-tA}= Q^*\begin{pmatrix} e^{-t\lambda_1}&z
{e^{-t\lambda_1}-e^{-t\lambda_2}\over \lambda_1-\lambda_2}\\ 0&
 e^{-t\lambda_2}\end{pmatrix}Q,
\label{2.11}
\end{equation*}
where
$$
{e^{-t\lambda_1}-e^{-t\lambda_2}\over \lambda_1-\lambda_2}=
-te^{-t\lambda_2}, \quad \text{if} \quad  \lambda_1 =\lambda_2.
$$
 Lemma \ref{l2.3} then yields
\begin{equation}
\begin{split}
\|e^{-tA}\|^2 &\leq\Tr(e^{-tA^*}e^{-tA})\\
&=e^{-2t Re(\lambda_1)}+e^{-2t Re(\lambda_2)}+|z|^2
\bigg|{e^{-t\lambda_1}-e^{-t\lambda_2}\over \lambda_1-\lambda_2}\bigg|^2.
\end{split}
\label{2.12}
\end{equation}
Therefore,  if $|\lambda_1-\lambda_2|>0$
\begin{equation*}
\|e^{-tA}\|^2 \leq (2+\frac{|z|^2}{|\lambda_1-\lambda_2|^2})e^{-2t
Re(\lambda_1)}
\end{equation*}
which proves part of the lemma.

Now, for  $|\lambda_1-\lambda_2|\geq\Lambda$, where
$\Lambda>0$ will  be chosen later,
\begin{equation}
{ |e^{-t\lambda_1}-e^{-t\lambda_2}| \over |\lambda_1-\lambda_2|}
\leq {2\over\Lambda} e^{-t Re(\lambda_1)};\label{useful}
\end{equation}
and for  $|\lambda_1-\lambda_2|\leq\Lambda$,
using $|e^\zeta-1|\leq |\zeta|\exp(|\zeta|)$ for any complex number
$\zeta$,
\begin{equation}
\begin{split}
|e^{-t\lambda_1}-e^{-t\lambda_2}|&=e^{-t \Re(\lambda_1)}
|e^{-t(\lambda_2-\lambda_1)}-1|\\
&\leq e^{-t\Re(\lambda_1)} t |\lambda_2-\lambda_1| e^{t\Lambda}.
\end{split}
\label{2.14}
\end{equation}
Therefore, choosing $\Lambda=\ds{1\over t}$ in
\re{useful} and \re{2.14},
\begin{equation*}
{ |e^{-t\lambda_1}-e^{-t\lambda_2}| \over |\lambda_1-\lambda_2|}
\leq Ct e^{-t Re(\lambda_1)}.
\end{equation*}
Substituting above into \re{2.12} completes the proof.
\endproof

\section{Decay rate of linear systems}\label{sec3}

In subsections \ref{sec3.1}, \ref{sec3.2} and \ref{sec3.3}, low-frequency 
($|\xi|$ close to $0$),  high-frequency (large $|\xi|$) and
middle range frequency analysis for the linear systems  are
performed respectively. 
We will identify  systems for which
there exist positive constants $\beta$ and $\delta_m$ such that
for any $t>0$ and
\begin{equation}
\begin{split}
\bullet & \text{ for } |\xi|\leq\delta_m,\ \|\exp(-tA)\|\leq C
\exp(-\beta t\xi^2), \\
\bullet & \text{ for }|\xi|>\delta_m,\  \|\exp(-tA)\|\leq
C\exp(-\beta t).
\end{split}\label{dichotomy}
\end{equation}
Here $\|\exp(-tA)\|$ is the norm of the linear operator
$\exp(-tA(\xi)) $ acting on $\mathbb{C}^2$. 
The generic constants $C$ and $\beta$ are independent of $t$ and $\xi$.
If  $\delta_m=+\infty$ is feasible in \re{dichotomy}, the system
is  in the {\it KdV--Burgers class}. Otherwise the system is  in the
{\it BBM--Burgers class}.
A summary of decay rates for the linear
systems is  given in subsection \ref{sec3.4}.

\begin{rmk} It is worth to note that \re{dichotomy} is sufficient but
  not necessary for proving the desired decay rate. An example is
 given in  Section \ref{sec3.3}
 where \re{dichotomy} is not valid but the linear system 
has the desired  decay rate.
\end{rmk}

\subsection{ Low frequency analysis}\label{sec3.1}

We now prove that for $|\xi|\to 0$, all systems are equivalent.  This is to say

\begin{pro}\label{pro3.1} There exists positive constants $\delta_m, \beta$ and $C$
depending   on the data
  $a,b,c,d$ and $\nu$,
such that for $|\xi|\leq\delta_m$ and for any $t>0$,
\begin{equation}
\|\exp(-tA)\|\leq C\exp(-\beta\xi^2 t ).\label{2.18}
\end{equation}
Consequently, for any initial data $Y_0$ with
$\text{supp}(\what  Y_0)\subset [-\delta_m,\delta_m]$,
\begin{equation*}
E(t)\leq C\ t^{-1/2} \|Y_0\|^2_{L^1_x}.
\label{2.19}
\end{equation*}
\end{pro}

\proof By referring to the  definitions of $\Delta$, $\alpha$ and $\ep$, one sees that
as $|\xi| \rightarrow 0$, 
$$\Delta\sim-4\xi^2 \quad \text{and} \quad \Tr(A)=\nu\alpha+\varepsilon\sim (\nu+1)\xi^2.$$
Therefore,  there exists $\delta_m > 0$ such that for $\xi$ in
 $[-\delta_m,\delta_m]$, $\Delta \leq 0$ and
$$ \half \leq \frac{\Tr(A)}{(\nu+1)\xi^2} \leq 2.$$
\re{2.18} then 
follows promptly from \re{2.15}.

Now, for any initial data $Y_0$ with
$\text{supp}(\what  Y_0)\subset [-\delta_m,\delta_m]$,
\begin{equation*}
\begin{split}
E(t)&=\int_{|\xi|\leq\delta_m } |\what Y (t,\xi)|^2d\xi
\leq C\int \exp(-2\beta t\xi^2)d\xi
(\sup_\xi(|\what{Y_0}(\xi)|^2))\\
&\leq C t^{-1/2}\| Y_0\|^2_{L_x^1},
\end{split}
\end{equation*}
by using the change of variable $\tau=\sqrt{2 \beta t} \xi$.
\endproof

\subsection{ High frequency analysis}\label{sec3.2}

The complete dissipation  and
the partial dissipation cases have to be studied separately. In the latter case, we will give one
example where the decay rate can be arbitrarily small.

Introducing the number
\begin{equation*}
\{r\}=\begin{cases}  1, & \quad \text{if } r\not= 0, \\
0, &\quad \text{if } r=0, \end{cases}
\end{equation*}
for $r \in \R$. Then
$\text{order}(\sigma)=\{a\}+\{c\}-\{b\}-\{d\}$
and $\text{order}(\ep)=2-2\{d\}$.

\subsubsection{ The complete dissipation case ($\nu=1$).}

It is observed in the
following that order$(\sigma)$ dictates if the system is in the
KdV-Burgers class or in the BBM-Burgers class.

\begin{pro} Assume $\nu=1$. For any $\delta>0$, there exists
$\beta >0$, such that if  $\text{supp}(\what Y_0)\subset \R \backslash
[-\delta,\delta]$,
\begin{equation}
E(t)\leq \exp(-2\beta t) \|Y_0\|^2_{L^2_x},\label{2.22}
\end{equation}
for any $t>0$.
In addition,
\begin{itemize}
\item if  $\order(\sigma)\leq 0$, the system is  in the {\it
BBM--Burgers class}. Namely, there exist positive constants
$\delta_M, \beta$ and $C$, such that for $|\xi| > \delta_M$ and
$t>0$,
$$ \| \exp(-tA)\| \leq C  e^{-\beta t};$$
\item if $\order(\sigma)\geq 1$,  the system is  in the {\it
KdV--Burgers class}. Namely, there exist positive constants
$\delta_M, \beta$ and $C$, such that for $|\xi| > \delta_M$ and
$t>0$,
$$\| \exp(-tA) \|\leq C e^{-\beta \xi^2 t}.$$
\end{itemize}
\label{pro1}\end{pro}

\proof
>From \re{2.3}, one finds that for  $\xi$ almost everywhere,
\begin{equation*}
{1\over 2} {d\over dt}  |\what Y(t,\xi)|^2+\alpha(\xi) |\what \eta(t,\xi)|^2
+\varepsilon(\xi) |\what w(t,\xi)|^2=0.
\label{2.21}
\end{equation*}
This gives directly, by setting $\beta= \min \{ \alpha(\delta),
\ep(\delta)\}$ which is positive, that \re{2.22} is valid. Furthermore, 
$\|\exp(-tA)\| \leq C  e^{-\beta t}$ for $|\xi| > \delta$. 
To figure out if the system is in
the BBM-Burgers or in the KdV-Burgers class, we separate the cases
as follows.

\begin{itemize}
 \item Assume first $\order(\sigma)\geq 1$. Then either $d=0$ or $b=0$.
Without loss of generality, let us assume $d=0$.

\begin{itemize}
\item{If $\Delta= (\alpha-\ep)^2-4\xi^2\sigma^2 >0$ for
$|\xi|$ large enough,} then  $\order(\sigma)=1$. In that case,
there exist $\beta >0$ and $\delta_M >0$
\begin{equation*}
\lambda_1\geq \frac{\det(A)}{\Tr(A)}=
\frac{\alpha+\sigma^2}{\frac {\alpha} {\xi^2} +1}\geq 2\beta \xi^2
\end{equation*}
for $|\xi| > \delta_M$. 
 By using  \re{2.16}
\begin{equation*}
\|\exp(-tA)\| \leq C (1+t\xi^2) e^{-2\beta t\xi^2}\leq  C
e^{-{\beta} t\xi^2}
\end{equation*}
for $|\xi| > \delta_M$ and the system is  in the KdV-Burgers class.
 \item{If $\Delta\leq 0$ for
$|\xi|$ large enough,} then there exists $\delta_M >0$ such that for $|\xi| > \delta_M$, 
$\half \xi^2\leq \Tr(A) \leq 2\xi^2$, and 
\re{2.15} implies  that the system is  in the KdV-Burgers class.
\end{itemize}

\item Assume now that $\order(\sigma)\leq 0$.
\begin{itemize}
\item{If $b\neq 0$ and $d\neq 0$ (weakly dispersive systems)} then
for $|\xi|$ large enough, $ {\rm Re}(\lambda_1)\leq \Tr(A) \sim
\frac{1}{b}+\frac{1}{d}$ as $|\xi| \rightarrow \infty$. This shows that a damping like
$e^{-\beta t\xi^2}$ is unlikely for high frequencies. Therefore  the
weakly dispersive systems are in the BBM-Burgers class.
\item{if
($b\neq 0$ and $d=0$) or ($b=0$ and $d\neq 0$)}. Without loss of
  generality, let us consider the case
$b\neq 0$ and $d=0$. Since $\Delta\sim \xi^4$ as  $|\xi| \to \infty$, we
have
\begin{equation*}
\lambda_1\leq \frac{2\det(A)}{\Tr(A)}=\frac {2(\alpha+\sigma^2)}{\frac {\alpha} {\xi^2} +1}
\sim C=O(1)
\end{equation*}
as  $|\xi| \to \infty$. This shows that a damping like $e^{-\beta t\xi^2}$ is again unlikely
for high frequencies. Therefore the system is in the BBM-Burgers class.
 \end{itemize}
\end{itemize}
\endproof

\subsubsection{ The partial dissipation case ($\nu=0$).}

We first note that when a system satisfies (C2) hypothesis, $\sigma=0$ and therefore 
$\lambda_1=0$ at $\xi=a^{-\half}$. But  one can always chose 
${\delta_M}$ large enough so  for $|\xi| > {\delta_M}$,  $\sigma$  is positive,
bounded from  below and  away from zero. Therefore the point where $\sigma$ vanishes will be 
considered in the next subsection.  We
now prove that in the partial
dissipation case, the decay rate is related to 
$\order(\sigma)$ and the 
strength of the dissipation which is characterized by  $\order(\ep)$.

\begin{pro} With  $\nu=0$, 
\begin{itemize}
\item if order$(\sigma)\geq 2-\frac{1}{2}$order$(\ep)$, then the
system is  in the KdV-Burgers class; 
\item if $|$order$(\sigma)|<
2-\frac{1}{2}$order$(\ep)$, then the system is  in the BBM-Burgers
class.
 \newline
In above two  cases,
when $\what Y_0$ is supported in $\Bbb R\backslash [-{\delta_M},{\delta_M}]$, then
for any $t>0$,
\begin{equation*}
E(t)\leq C\exp(-2\beta t)\|Y_0\|^2_{L^2_x}; \label{2.23}
\end{equation*}
\item if order$(\sigma)\leq -2+\frac{1}{2}$order$(\ep)$,
arbitrarily slow decay can occur. \label{pro2}
\end{itemize}
\end{pro}

\proof
To begin, one observes that    $\Delta=\ep^2-4\xi^2\sigma^2$. When 
$\Delta >0$, $\lambda_1$ satisfies
\begin{equation}\label{directestimate}
2\lambda_1=\varepsilon \left(1-\left(1-
{4\xi^2\sigma^2\over\varepsilon^2}\right)^{1/2}\right)
\end{equation}
which is a direct consequence of \re{2.7b}.

\begin{itemize}
\item {\it When  order$(\ep)$=0, i.e $d\neq 0$ (and
$\order(\sigma)\leq 1$):}
\begin{itemize}
\item if $\order(\sigma)\geq -1$, and  if
$\Delta=\varepsilon^2-4\xi^2\sigma^2>0$ for $|\xi|$ large enough, which  is
possible only for $\order(\sigma)=-1$, we have
\begin{equation*}
\lambda_1\geq \frac{\det(A)}{\Tr(A)}\geq C\xi^2\sigma^2=O(1)
\end{equation*}
as $|\xi| \to \infty$.
Then by \re{2.16} we have
\begin{equation*}
\|e^{-tA}\|\leq C(1+t|\xi|\sigma)e^{-t\lambda_1}\leq Ce^{-\beta t}
\end{equation*}
for $|\xi|$ large enough and the system is in the BBM-Burgers class. On the other hand, if
$\Delta\leq 0$ for high frequencies, since $\Tr(A)\sim
\frac{1}{d}$ as $|\xi| \rightarrow \infty$, then \re{2.15} implies that the system is  in the
BBM-Burgers class; 
\item if $\order(\sigma)=-2$, $\Delta \sim
\frac{1}{d^2} >0$ as  $|\xi| \rightarrow \infty$ and  by \re{directestimate}, $\lambda_1\sim
C|\xi|^{-2}$, 
  therefore arbitrarily  slow decay could  occur. An
example of such case will be given below.
\end{itemize}
\item {\it When $\order(\ep)=2$ i.e $d=0$ (and $\order(\sigma)\geq
-1$):} $\Delta=\xi^2(\xi^2-4\sigma^2)$ has a limit $\Delta_0$ in
$[-\infty, +\infty]$ when $|\xi|$
approaches  $+\infty$.
\begin{itemize}
\item If $\Delta_0$  is in $(-\infty, 0]$, then since $\Tr(A)=\xi^2$,  \re{2.15} implies
that the system is in the KdV-Burgers class. This occurs when 
$\order(\sigma)=2$ and may occur when  $\order(\sigma)=1$.

 \item If $\Delta_0$ is in $(0, +\infty]$, then
since
\begin{equation*}
\frac{2\det(A)}{\Tr(A)}\geq \lambda_1 \geq
\frac{\det(A)}{\Tr(A)} =\sigma^2, 
\end{equation*}
\re{2.16} implies for any $\xi$ 
\begin{equation}\label{3.15}
\|\exp (-tA)\|\leq C\left(1+|\xi| \sigma\min\left(t,
\frac{1}{\sqrt{\Delta}}\right)\right)\exp(-t \sigma^2).
\end{equation}
\begin{itemize}
\item If $\order(\sigma)=1$,  then \re{3.15} implies the system is  in the KdV-Burgers class.
 \item If $\order(\sigma)=0$, 
$\frac{|\xi\sigma|}{\sqrt{\Delta}}=0(1)$ as $|\xi| \rightarrow \infty$, the system is
in the BBM-Burgers class. And similarly, 
 \item if $\order(\sigma)=-1$,  any arbitrarily
slow decay could  occur.

\end{itemize}
\end{itemize}
\end{itemize}
\endproof

{\bf Example of slow decay:} Consider the linearized  BBM-BBM system with
partial dissipation,
\begin{equation*}
\begin{split}
& \eta_t+u_x-b \eta_{xxt}=0,\\
& u_t+\eta_x-d u_{xxt}=u_{xx},
\end{split}
\label{2.25}
\end{equation*}
which has $\order(\ep)=0$ and $\order(\sigma)=-2$.
Since as  $|\xi|\to+\infty$,
\begin{equation*}
\begin{split}
\Delta& \sim{1\over d}>0, \quad 
|z|=\ds{2|\xi|\sigma\sim{2\over\sqrt{bd}}\ {1\over |\xi|}},\\
2\lambda_1&=\Tr(A)\left(1-\left(1-{4\det(A)\over (\Tr(A))^2}\right)^{1/2}\right)
\sim\ 2\ {\det(A)\over \Tr(A)}\sim {2\over d\xi^2}.
\end{split}
\end{equation*}
 Therefore
\begin{equation*}
\ds \|e^{-tA}\|\leq C\exp \left(-{\beta t\over\xi^2}\right) \label{2.26}
\end{equation*}
which shows that any arbitrary
slow decay could occur.

\subsection{Middle range frequency analysis}\label{sec3.3}

We first note from Lemma \ref{lem2.6} that to get the optimal decay estimate for 
 the cases where $\det(A)$ (and therefore $\lambda_1$)
has a zero for $|\xi|>0$, these cases  need to be discussed separately.
Therefore, we have the following two propositions. 

\begin{pro}\label{pro3.2}
Assume that $\nu=1$, or that $\nu=0$ and the dispersive coefficients $a, b, c, d$
satisfy (C1). Then for any $\delta_m$ and $\delta_M$,  $0<\delta_m \leq \delta_M$,
 there exists $\beta>0$ such that for
$|\xi| \in [\delta_m,\delta_M]$ and for any $t>0$,
\begin{equation}
\|\exp (-t A)\| \leq C \exp (-\beta t). \label{middle}
\end{equation}
Moreover for any $\what{Y_0}$ with support included in
$[\delta_m,\delta_M]\cup [-\delta_M,-\delta_m]$,
$$E(t)\leq C \exp (-2\beta
t)||Y_0||^2_{L_x^2}.$$
\end{pro}

\proof  Since  $\Tr{(A)}$ and $\det(A)$  cannot vanish for $|\xi|
\in [\delta_m,\delta_M]$
 under the assumptions,
\re{middle} is the direct consequence of
\re{2.15} and \re{2.16}. In addition
\begin{equation*}
E(t)\leq \sup_\xi ||e^{-tA}||^2 \, ||\what{Y_0}||^2_{L^2_\xi}\leq C
\exp (-2\beta t) ||Y_0||^2_{L_x^2},
\end{equation*}
which completes the proof of the proposition.
\endproof

\begin{rmk}\label{regledetrois}
By noticing  that \re{middle} can be replaced by
\begin{equation*}
\|\exp (-t A)\| \leq C \exp (-\beta^* t\xi^2). \label{middleb}
\end{equation*}
with $\beta^*=\beta/\delta_M^2$,
the middle range frequency analysis and  the high frequency analysis
can be combined to simplify
certain calculations  for these systems regardless if they are in BBM-Burgers class or KdV-Burgers class.
\end{rmk}

\begin{pro}\label{pro3.6} Assume that $\nu=0$ and the dispersive coefficients satisfy  (C2).
  Then for any $0<\delta_m<\delta_M$ with  $r=a^{-\half} \in [\delta_m, \delta_M]$,  there
exists $\beta>0$ and $C>0$ such that for any $|\xi| \in
[\delta_m,\delta_M]$ and for any $t>0$,
\begin{equation}\label{badcase}
\|\exp (-t A)\| \leq C {\rm exp}\{-\beta t(|\xi|-r)^2) \}.
\end{equation}
Moreover for any $\what{Y_0}$ with support included in
$[\delta_m,\delta_M]\cup [-\delta_M,-\delta_m]$ and for any $t >0$,
$$E(t)\leq Ct^{-1/2}||Y_0||^2_{L_x^1} .$$
\end{pro}

\begin{rmk}
Proposition \ref{pro3.6}  shows that even when the dichotomy is not valid,
the energy could decay as O$(t^{-1/4})$ when $t$ goes to $+\infty$.
\end{rmk}
\proof Since  $\det (A)$ vanishes at $r=(\sqrt{a})^{-1}$ and, when
$|\xi| \to r$, $\Delta\sim \frac{a}{a+d}>0$,  $\lambda_1\sim \beta {\rm
det}(A)\sim \beta \sigma^2 \sim \beta (|\xi|-r)^2$. 
 Therefore, from \re{2.16},
\begin{equation*}
\|\exp (-t A)\| \leq C (1+ \min(t,1)\sigma) {\rm
exp}(-\beta \sigma^2t) \label{3.7b}
\end{equation*}
as $|\xi|$ in the neighborhood of $r$. 
Using the fact that for $t>0$, $\min(t,1)\leq \sqrt{t}$, so 
there exists $C>0$
such that 
$$\|\exp (-t A)\| \leq C(1+\sqrt{t}\sigma){\rm exp}(-\beta \sigma^2t)\leq C{\rm
exp}(-\frac\beta{2} \sigma^2t),$$
 we obtain the estimate \re{badcase} for $|\xi|$ close to
$r$. For other $|\xi|$ in $[\delta_m, \delta_M]$,  the same argument in  the proof
of Proposition \ref{pro3.2} and Remark \ref{regledetrois} applies. For the decay rate of
$E(t)$, same argument as in the proof of Propositions
\ref{pro3.1} can be used. In fact, $|\xi|-r$ plays the same role as $|\xi|$ in that
case.
\endproof

\subsection{Decay for linear systems}\label{sec3.4}

Since linear system \re{1.6} defines a  semi-group 
$e^{-tA}$ for  $t\geq 0$, that is  contracting on $L^2\times L^2$ in 
the variable $(\eta, w)$, 
the initial value problem  is therefore well-posed and
the $L^2$ norm decays. 

Combining the low, middle and high
frequency analysis, the decay rate for the linear system \re{1.6}
can be stated as
\begin{thm}\label{thm3.8} For systems \re{1.6} with  the dispersive 
constants $a,b,c,d$ satisfy the constraints
(C0)-(C1) or (C0)-(C2), assuming either $\{\nu=1\}$ or,
$\{\nu=0$ and
$\order(\sigma)>-2+\frac{1}{2}\order(\varepsilon)\}$, then for any
$(\eta_0, Hu_0)=(\eta_0, w_0) \in (L^1(\R) \cap L^2(\R))^2$ where
$\what H \what u_0= (\frac
{(1-a\xi^2)(1+d\xi^2)}{(1-c\xi^2)(1+b\xi^2)})^\half \what u_0$,
there exists a constant $C$, such that for any $t>0$
$$ \|(\eta, Hu)\|_{L_x^2}=\|(\eta,
w)\|_{L_x^2} \leq C t^{-1/4}.
$$
\end{thm}

\begin{rmk}
 This is equivalent to say ,
with respect to physical variables $(\eta, u)$, that for any
$(\eta_0, u_0) \in (L_x^2\cap L^1_x) \times (H_x^h \cap
W_x^{h,1})$, 
$$
\|(\eta, u)\|_{L_x^2 \times H_x^h} \leq C t^{-1/4}
$$
for any  $t>0$, where $h=\order(\what{H})=\{a\}+\{d\}-\{c\}-\{b\}$.
\end{rmk}

\proof
Combining the low, middle and high frequency analysis, we have
$$E(t)=E_{low }(t)+E_{middle}(t)+E_{high}(t)
\leq C(\eta_0, u_0)(t^{-\half}+ e^{-2\beta t})$$ 
for $t>0$  
 where  $C(\eta_0, u_0)$ is a function of the dispersive 
coefficients and the norms of $\eta_0$ and $u_0$.
\endproof

We complete this section by the following

\begin{cor} For any dissipation, the classical Boussinesq
 system,  the Bona--Smith system, the
 coupled KdV--BBM ($b=c=0$) system, 
the BBM--KdV systems ($a=d=0$) and the 
weakly dispersive systems ($b>0$ and
$d>0$) with  $a<0$ or $c<0$
 belong to the {\it BBM--Burgers
class}.
\end{cor}

\begin{cor} With complete dissipation,
the  KdV--KdV system ($b=d=0$, $a=c>0$) belongs to the {\it
KdV--Burgers class};  the weakly dispersive systems ($b>0$ and
$d>0$) belong to the {\it BBM--Burgers class}.
\end{cor}

\begin{rmk} When the consideration is restricted to the linear systems,
a result  similar to Theorem \ref{thm3.8} (substituting  $\alpha$ for $\ep$ in the
statement)  holds
when one replaces
 $(\nu \eta_{xx},u_{xx})$ by $(\eta_{xx},\nu u_{xx})$
 in the right-hand side of  \re{2.1} since
 $\eta$
and $u$ play the same role.
\end{rmk}

\section{Nonlinear Theory}\label{sec4}

For convenience, we will only consider in this section (i) the
complete dissipation and (ii) the partial dissipation with $a, b,
c, d$ satisfy the (C1) assumption. The partial dissipation with
$a,b,c,d$ satisfy (C2) will be considered elsewhere.

\subsection{ A general result}

Consider an evolution equation that reads
\begin{equation}
v_t+Lv+\partial_x(F(v))=0, \quad v(t=0)=v_0\label{3.1}
\end{equation}
where $L$ is a  linear unbounded operator with symbol $A$  
and $F$ is a nonlinear
quadratic operator.

Assuming that $L$ generates a semi--group $S(t)$ on $L^2_x$ that
satisfies the {\it dichotomy assumption} (\ref{dichotomy}), namely
there exist $\beta>0$ and $\delta>0$  ($\delta$ can be $+\infty$) such that
for any $t>0$ and 
\begin{equation}
\begin{split}
\bullet & \text{ for } |\xi|\leq\delta,\ \|S(t)\|\leq C
\exp(-\beta t\xi^2),\\
\bullet & \text{ for }|\xi|>\delta,\  \|S(t)\|\leq C  \exp(-\beta
t).
\end{split}\label{dich2}
\end{equation}
In addition,
\begin{equation}
\sup_{t\geq 0}(t^{1/4}\|S(t) v_0\|_{L^2_x})\leq C_1 \|v_0\|_{L_x^1\cap L_x^2}=\overline
C.\label{3.15b}
\end{equation}
Assuming  also  the nonlinear term  $F(v)$ satisfies
\begin{equation}
\sup_{|\xi|\leq\delta} |\what  F|+\left(\int_{|\xi|\geq\delta} |\xi|^2
|\what  F|^2 d\xi\right)^{1/2}
\leq C\|v\|^2_{L^2_x}
\label{3.4}
\end{equation}
for any $t>0$, where $\what F= \mathcal{F}(F(v))$.

 Let us recall that a mild
solution to
\re{3.1} is a solution to the integral equation
\begin{equation}
v(t)=S(t)v_0-\int_0^t S(t-s)\partial_x F(v(s))ds.\label{3.7}
\end{equation}
Under the above assumptions, we may construct a solution to
  \re{3.7}
 by performing a fixed point argument (see \cite{CH98}, \cite{Karch99},
\cite{BCS04}, \cite{CazWei00}) on 
the space
$$ E=\left\{ u: \ \sup_{t>0} \{t^{1/4}\|u(t)\|_{L^²_x}  \} 
< \infty\right\}, $$
which is a Banach space of functions that are continuous
in time with value in $L^2$ that are $O(t^{-1/4})$
when $t$ goes to $+ \infty$.
If $\overline C$  is small enough, 
a fixed point argument to the Duhamel's form
of the equation in the ball  in $E$ centered at origin  would 
provide the solution.

\begin{thm}\label{thm4.1} 
 For system \re{3.1} with assumptions
\re{dich2}-\re{3.15b}-\re{3.4}, 
there exists a numerical constant $ C$ such that
  for
 any {\it mild solution} to \re{3.1} starting from $v_0$ with 
\begin{equation*}
 \|v_0\|_{L_x^1\cap L^2_x}\leq  C,\label{3.5}
\end{equation*}
then
\begin{equation}
\|v(t)\|_{L^2_x} \leq O(t^{-1/4}) \quad \text{as} \quad t \rightarrow \infty.\label{3.6}
\end{equation}
\end{thm}

\proof
To begin, we first control  the low frequency part of the nonlinear term.
Let
\begin{equation*}
\begin{split}
\what  N&:={\mathcal F}\left(\int_0^t S(t-s)\partial_x F(v(s))ds\right)\\
&=i\int_0^t e^{-(t-s)A}\xi\what  F(v(s))ds.\label{3.8}
\end{split}
\end{equation*}
Using the first inequality in  \re{dich2} in combination with
\re{3.4}, one obtains 
\begin{equation}
\begin{split}
& \left(\int_{|\xi|\leq\delta} |\what  N|^2 d\xi\right)^{1/2}\leq
C \int_0^t\left[\int_{|\xi|\leq\delta}\|e^{-(t-s)A}\|^2\xi^2|\what  F|^2
 d\xi\right]^{1/2}ds\\
& \hskip .5cm  \leq C\int_0^t\left[\int_{\R} \xi^2
e^{-2\beta(t-s)\xi^2}d\xi\right]^{1/2}
\left(\sup_{|\xi|\leq\delta}|\what  F|\right) ds\\
& \hskip .5cm \leq C\int_0^t {\|v(s)\|^2_{L^2_x}\over (t-s)^{3/4} }\ ds.
\end{split}
\label{3.9}
\end{equation}
We now control  the high frequency part of the nonlinear term,
using the second inequality in \re{dich2} in combination with
\re{3.4}
\begin{equation}
\begin{split}
\left( \int_{|\xi|\geq\delta} |\what  N|^2 d\xi\right)^{1/2}&\leq
C\int_0^t
e^{-\beta(t-s)}\left( \int_{|\xi|>\delta}\xi^2 |\what  F|^2 d\xi\right)^{1/2} ds\\
&\leq C\int_0^t e^{-\beta(t-s)} \|v(s)\|^2_{L^2_x} ds.
\end{split}
\label{3.10}
\end{equation}
Introducing  the norm
\begin{equation}
M(t)=\sup_{s\in [0,t] } (s^{1/4} \|v(s)\|_{L^2_x}), \label{3.11}
\end{equation}
and if $v$ solves \re{3.7}, then due to \re{3.9}--\re{3.10}, 
\begin{equation*}
\begin{split}
& t^{1/4} \|v(t)\|_{L^2_x}\leq t^{1/4} \|S(t) v_0\|_{L^2_x}+t^{1/4} \| \hat N(t) \|_{L_{\xi}^2}\\
&\leq  t^{1/4} \|S(t) v_0\|_{L^2_x}
 + CM(t)^2\int_0^t\left[{t^{1/4}\over s^{1/2}(t-s)^{3/4}}+{t^{1/4}\over
      s^{1/2} }
\exp (-\beta(t-s))\right] ds.
\end{split}
\label{3.12}
\end{equation*}
By applying the change of variable $s=\tau t$ in the   integration, one
finds
\begin{equation*}
\begin{split}
&\int_0^t\left[{t^{1/4}\over s^{1/2}(t-s)^{3/4}}+{t^{1/4}\over
      s^{1/2} }
\exp (-\beta(t-s))\right] ds
\\
\leq  & C+ \int_0^1 \frac{t^{3/4}}{\tau^\half} \exp(-\beta
t(1-\tau)) d\tau \leq C.
\end{split}
\end{equation*}
Therefore, using the
property  \re{3.15b},
the positive, nondecreasing function $M(t)$ satisfies $M(0)=0$ and for any
$t\geq 0$,
\begin{equation}
C_0 M(t)^2-M(t)+\overline C\geq 0, \label{3.16}
\end{equation}
where $C_0$ is a positive constant. Choosing $\overline C$ such that
$C_0 x^2-x+\overline C= 0$ has two real roots $0<r_1<r_2$, namely choosing $\overline C < \frac 1{4C_0}$,
then \re{3.16} holds only if
$M(t)$ is trapped in the interval $[0,r_1]$.
Therefore when
$$\|v_0\|_{L_x^1\cap L_x^2} \leq \frac {\overline C}{C_1} < \frac 1{4C_0C_1}$$
$M(t)$ is bounded and \re{3.6} is valid. \endproof

\subsection{  Applications to weakly dispersive systems with complete dissipation
or partial dissipation}

Since  $b>0$ and $d>0$, $\order(\sigma) \leq 0$ and the corresponding 
linearized system is in the BBM-Burgers class.
From
Proposition \ref{pro1} and Proposition \ref{pro2}, this corresponds to consider {\it
complete
  dissipation},
 or {\it partial dissipation} together with $\order(\sigma)\geq -1$.

\begin{thm}\label{thm4.2} Consider a weakly dispersive two--way wave model
  $(b>0 \text{ and } d>0)$
with either the  complete dissipation or the  partial dissipation
together with $a<0$ or $c<0$. Then, for small initial data, 

\begin{itemize}
\item if $H$ is of order 0,
\begin{equation}
\|\eta(t)\|^2_{L^2_x}+\|u(t)\|^2_{L^2_x} \leq O(t^{-1/2});\label{3.26}
\end{equation}
\item if $H$ is of order 1,
\begin{equation}
\|\eta(t)\|^2_{L^2_x}+\|u(t)\|^2_{H^1_x} \leq O(t^{-1/2});\label{3.27}
\end{equation}
\item if  $H$ is of order $-1$,
\begin{equation}
\|\eta(t)\|^2_{H^1_x}+\|u(t)\|^2_{L^2_x} \leq O(t^{-1/2});\label{3.28}
\end{equation}
\end{itemize}
as $t \to \infty$.
\end{thm}

\proof
Note that the theorem is proved after \re{3.4} is validated and we will do that 
by 
discussing  the cases according to the order of $H$. 
\begin{itemize}
\item If  $H$ is of order 0 or 1.
Introducing the change of variable
\begin{equation*}
v=\mathcal{F}^{-1}(\what \eta, \what H\what  u)= \mathcal{F}^{-1}(\what
\eta, \what w),\label{3.17}
\end{equation*}
the full nonlinear system
\begin{equation}
\begin{split}
& \eta_t+u_x+au_{xxx}-b\eta_{xxt}+(\eta u)_x=\nu\eta_{xx},\\
&u_t+\eta_x+c\eta_{xxx}-du_{xxt}+uu_x=u_{xx},
\end{split}
\label{3.18}
\end{equation}
transforms to
$$
v_t+Lv=-\partial_x F(v),
$$
where  $L$ has  symbol  $A=\begin{pmatrix} \nu\alpha&i{\rm sgn}(\omega_1)\xi\sigma\\
i{\rm sgn}(\omega_2)\xi\sigma&\varepsilon\end{pmatrix}$ and
$F$ reads
\begin{equation*}
F(v)=\begin{pmatrix} (1-b\partial^2_x)^{-1}\eta H^{-1}w\\
\frac{1}{2}H(1-d\partial^2_x)^{-1}(H^{-1}w)^2\end{pmatrix}.\label{3.19}
\end{equation*}
To check \re{3.4}, it is natural  to separate the estimate into two parts. 

\begin{itemize}
\item {\it Low frequency ($|\xi| \leq \delta$) estimate}:
Since $H^{-1}$, which has order $0$ or
$-1$, is bounded on $L^2_x$, straightforward computations lead to
\begin{equation*}
\begin{split}
\|\what{\eta H^{-1}
w}\|_{L^\infty_\xi}&+\|\what{(H^{-1}w)^2}\|_{L^\infty_\xi}
\leq C(\|\eta\|^2_{L^2_x}+\|H^{-1}w\|^2_{L^2_x})\\
&\leq C(\|\eta\|^2_{L^2_x}+\|w\|^2_{L^2_x}). \label{3.20}
\end{split}
\end{equation*}
Because
 $(1-b\partial_x^2)^{-1}$ and
$H(1-d\partial_x^2)^{-1}$
 are bounded operators,
$$\sup_{|\xi|\leq\delta}|\what  F| \leq C \|v\|^2_{L_x^2}.$$

\item {\it High frequency ($|\xi| > \delta$) estimate}: Since
$\partial_x(1-b\partial^2_x)^{-1}$ is a smoothing operator
(it is of order $-1$) and $H^{-1}$ is a bounded operator on $L^2_x$, 
\begin{equation*}
\begin{split}
& \|\partial_x (1-b\partial^2_x)^{-1}(\eta H^{-1}w)\|_{L^2_x} \leq 
C\|\eta H^{-1}w\|_{H^{-1}_x} \\
& \hskip .5cm \leq  C\|\eta\|_{L^2_x} \|H^{-1}w\|_{L^2_x}
\leq C(\|\eta\|^2_{L^2_x}+\|w\|^2_{L^2_x}),
\end{split}
\end{equation*}
where Lemma 2.2(ii) in \cite{BCS04},
$$ \| fg \|_{H^{-1}} \leq C \|f\|_{L^2} \|g\|_{L^2} $$
is used. 
Now, consider  the term 
$\partial_x(H(1-d\partial^2_x)^{-1})((H^{-1}w)^2)$ 
in $F$.  If $H$ is of order $0$, it can be
 bounded in the same way.
If $H$ is of order $1$, then
\begin{equation*}
\begin{split}
& \|\partial_x H(1-d\partial_x^2)^{-1} (H^{-1}w)^2\|_{L^2_x} \leq 
C\|(H^{-1}w)^2\|_{L^2_x} \\
& \hskip .5cm \leq  C\|H^{-1}w\|^2_{H^1_x}\leq C\|w\|^2_{L^2_x},\label{3.22}
\end{split}
\end{equation*}
where Lemma 2.2(iv) in \cite{BCS04}
$$ \| fg \|_{L^{0}} \leq C \|f\|_{H^1} \|g\|_{H^1} $$
is used. 
\end{itemize}

Combining the lower frequency and higher frequency analysis, one sees \re{3.4} is valid and Theorem
\ref{thm4.1} yields the desired result.

\item  If $H$ is of order $-1$, introducing  the change of variable
\begin{equation}
v=\mathcal{F}^{-1}((H^{-1}\what \eta,\what  u)).\label{3.23}
\end{equation}
and setting  $\what  \tau=H^{-1}\what \eta$,  the full nonlinear
system \re{3.18} reads
\begin{equation*}
v_t+Lv=-\partial_x F(v)\label{3.24}
\end{equation*}
where $L$ has symbol $A=\begin{pmatrix} \nu\alpha&i{\rm sgn}(\omega_1)\xi\sigma\\
 i{\rm sgn}(\omega_2)\xi\sigma&\varepsilon\end{pmatrix}$ and 
\begin{equation*}
F(v)=\begin{pmatrix} H^{-1}(1-b\partial^2 x)^{-1}(uH\tau)\\
\frac{1}{2}(1-d\partial^2_x)^{-1}(u^2)\end{pmatrix}.\label{3.25}
\end{equation*}
The proof is then very similar to the previous case  and therefore 
omitted.
\end{itemize}
\endproof

\begin{cor} For the following two special cases, we have 

\begin{itemize}
\item solutions to  Bona--Smith system 
($a=0, b>0, c<0$ and  $d>0$) with complete or partial dissipation satisfy \re{3.28};

\item solutions to  BBM--BBM system with complete dissipation
($a=c=0, \nu=1)$ satisfy \re{3.26}.
\end{itemize}
\end{cor}

\subsection{Application to KdV-Burgers systems with complete dissipation.}

Using Proposition \ref{pro1}, 
this implies that order$(\sigma) \geq 1$, and then that $b=0$ and/or $d=0$.

\smallskip
\noindent First case: 
Consider the case where $b=d=0$.
The analysis in \cite{BCS04} implies that $a=c=1/6$, so
the system satisfies (C2) assumptions.
Since the dichotomy assumption \re{dich2} was proved  
in Section \ref{sec3} and the linearized system is in the  
{\it
KdV-Burgers class}, Theorem \ref{thm4.1} applies when 
\re{3.4} with  $\delta=+\infty$ is verified.

Let us observe that $H=1$ and that the full nonlinear system reads
\begin{equation*}
v_t+Lv=-\partial_x F(v)\label{3.29}
\end{equation*}
where $v=(\eta,u)$, $L$ has symbol $A=\begin{pmatrix} \xi^2&i{\rm sgn}(\omega_1)\xi\sigma\\
i{\rm sgn}(\omega_2)\xi\sigma&\xi^2\end{pmatrix}$ and 
\begin{equation*}
F(v)=\begin{pmatrix}\eta u\\ {u^2\over
2}\end{pmatrix}.\label{3.30}
\end{equation*}
Since 
\begin{equation*}
\|F(v)\|_{L^1_x}\leq C(\|\eta\|^2_{L^2_x}+\|u\|^2_{L^2_x}), \label{3.31}
\end{equation*}
\re{3.4} with $\delta=\infty$ is a direct consequence. 
\smallskip

\noindent Second case: Consider the case where
$b>0$ and $d=0$. Due to (C0), $c \geq 0$. Since $\order(\sigma) \geq 1$, the system must have  $a=c>0$.
Therefore order$(H)=-1$. Introducing the change of variable
\re{3.23} the system reads as \re{3.29}
with the following nonlinearity

\begin{equation*}
F(v)=\begin{pmatrix} H^{-1}(1-b\partial^2 x)^{-1}(uH\tau)\\
\frac{u^2}{2}\end{pmatrix}.\label{new2}
\end{equation*}
 \re{3.4} with $\delta=\infty$
is valid  and the proof is straightforward and then omitted.

\smallskip

\noindent Third case: Consider the case where
$b=0$ and $d>0$. Since order$(\sigma)=1$, then
$a$ and $c$ can not vanish and order$(H)=1$. 
In this case, with the change of variable
\re{3.17}, the system reads as \re{3.29}
with $F$ being

\begin{equation*}
F(v)=\begin{pmatrix} \eta H^{-1}w\\
\frac{1}{2}H(1-d\partial^2_x)^{-1}(H^{-1}w)^2\end{pmatrix}.\label{new3}
\end{equation*}
\re{3.4} with $\delta=\infty$
is again valid and the proof is straightforward and then omitted.

Therefore we can state

\begin{thm}\label{thm4.4} Consider a 
KdV-Burgers system with complete dissipation.
Then for small initial data, 

\begin{itemize}
\item if $H$ is of order 0,
\begin{equation}
\|\eta(t)\|^2_{L^2_x}+\|u(t)\|^2_{L^2_x} \leq O(t^{-1/2});\label{new4}
\end{equation}
\item if $H$ is of order 1,
\begin{equation}
\|\eta(t)\|^2_{L^2_x}+\|u(t)\|^2_{H^1_x} \leq O(t^{-1/2});\label{new5}
\end{equation}
\item if  $H$ is of order $-1$,
\begin{equation}
\|\eta(t)\|^2_{H^1_x}+\|u(t)\|^2_{L^2_x} \leq O(t^{-1/2});\label{new6}
\end{equation}
\end{itemize}
as $t \to \infty$.
\end{thm}

\begin{rmk}This is not surprising for KdV-KdV system ($b=d=0$)
since if we introduce the new
variables $\eta=w^1+w^2$ and $u=w^1-w^2$, then
 \re{3.18}
 reads as a system of two linear KdV--Burgers systems (weakly) coupled
 through
nonlinear terms.
See Section 2.3 in \cite{BCS04}.
\end{rmk}

\subsection{Other cases}

In some other cases, the  method presented here  does not work
straightforwardly.  As pointed out in the Introduction, other methods exist
which might enable the analysis to go further. 
These methods could  also be helpful  to extend
our local results to global ones.

\section{The $L^\infty_x$--decay rate}\label{sec5}

First, we observe that for the cases of weakly dispersive wave equations and 
 KdV-KdV system, the nonlinear terms  satisfy
\begin{equation}
\sup_{|\xi|\leq\delta} |\xi\what
F|+\left(\int_{|\xi|>\delta}|\xi|^4 |\what  F|^2 d\xi\right)^{1/2}
\leq C\|v\|_{L^2_x}\|v_x\|_{L^2_x}. \label{3.33}
\end{equation}

We now estimate the decay rate of
$\partial_x v(t)$ in $L^2_x$ when $v$ solves \re{3.7}. To begin, 
we differentiate \re{3.7} with respect to $x$ and treat 
the nonlinear term of the resulting equation with a  procedure similar
to the one in the proof of 
Theorem \ref{thm4.1}, but  using \re{3.33} instead of \re{3.4}. 

We first note that for  low frequencies, 
\begin{equation*}
\left(\int_{|\xi|\leq\delta} \xi^2 |\what  N|^2 d\xi\right)^{1/2}
\leq C\int_0^t {\|v(s)\|_{L^2_x}\|v_x(s)\|_{L^2_x}\over
(t-s)^{3/4} } ds \label{3.34}
\end{equation*}
and for high frequencies
\begin{equation*}
\left(\int_{|\xi|\geq\delta}\xi^2 |\what  N|^2
d\xi\right)^{1/2}\leq C\int_0^t e^{-\beta(t-s)}\|v\|_{L^2_x}
\|v_x\|_{L^2_x} ds.\label{3.35}
\end{equation*}

We now consider the linear part.
For $v_0$ in $H^1(\Bbb R)\cap L^2(\Bbb R)$, the linear term can be estimated by splitting 
the region of integration into low
frequencies and high frequencies. Using \re{dich2},  we have
\begin{equation*}
\begin{split}
\|\partial_x S(t)v_0\|^2_{L^2_x} &\leq C\bigg[\left(\int\xi^2
e^{-\beta \xi^2 t}d\xi\right)
\|v_0\|^2_{L^2_x}
+e^{-2\beta t}
\|\partial_x v_0\|^2_{L^2_x}\bigg]\\
&\leq {C\over t^{3/2}} \|v_0\|^2_{L^2}+Ce^{-2\beta
t}\|v_0\|^2_{H^1}.
\end{split}
\label{3.36}
\end{equation*}
Therefore, the linear part behaves like $O(t^{-3/4})$  as $t\to \infty$. 
Since $\| v\|_{L^2}$ is $O(t^{-1/4})$ as $t \to \infty$, we have 
\begin{equation*}
\begin{split}
 \|\partial_x v(t)\|_{L^2_x}&\leq C(v_0) t^{-3/4} 
 +C\sup_{s\in[0,t]} \left(s^{1/4}\|v(s)\|_{L^2}\right) \\
&\times \int_0^t {\|v_x(s)\|_{L^2_x}\over (t-s)^{3/4}}\ {1\over
s^{1/4}}+{e^{-\beta(t-s)}\over s^{1/4}}
\|v_x(s)\|_{L^2_x}ds.
\end{split}
\label{3.37}
\end{equation*}
Simple calculations show 
\begin{equation*}
\int_0^t \ {ds\over s^{1/4}(t-s)^{3/4}} +\int_0^t
{e^{-\beta(t-s)}\over s^{1/4}}\ ds\leq C.\label{3.38}
\end{equation*}
Therefore
\begin{equation}
\begin{split}
\|\partial_x v(t)\|_{L^2_x}&\leq C(v_0) t^{-3/4}
+C\sup_{s \in [0,t]}\|\partial_x v(s)\|_{L^2_x}
\sup_{s \in [0,t]} (s^{1/4}\|v(s)\|_{L^2}) \\
&\leq C(v_0) t^{-3/4}+C_2 M(t)\sup_{s \in [0,t]}\|\partial_x v(s)\|_{L^2_x}  .
\end{split}
\label{3.39}
\end{equation}
Since $M(t)$, that is defined
in \re{3.11}, is bounded by the first root $r_1$  of
$$
C_0 x^2-x+\overline C=0,
$$
and 
$
r_1 \sim \overline C$ as  $\overline C \to 0$ 
(since ${1\over C_0}(1-(1-4 C_0\overline C)^{1/2})\sim 2\overline C$
as $ \overline C \to 0$).
Therefore, when $\overline C$ small enough, $M(t) \leq 2 \overline C$. Hence by choosing 
$\overline C$ such that 
\begin{equation}
2C_2\overline{C}\leq {1\over 2},\label{3.40}
\end{equation}
\re{3.39} leads to
\begin{equation}
\|\partial_x v(t)\|_{L^2_x}\leq 2 C(v_0) t^{-3/4}\label{3.41}
\end{equation}
and we obtain the following theorem. 

\begin{thm}\label{thm5.1}
 For system \re{3.1} with assumptions
\re{dich2}-\re{3.15b}-\re{3.4}, 
assume $v_0$ is in $H^1(\Bbb R)\cap L^1(\Bbb R)$ and
$\|v_0\|_{L^1\cap L^2}$ is small enough.
Then
\begin{equation*}
\|v\|_{L^\infty_x}\leq O(t^{-1/2})\label{3.42}
\end{equation*}
as $t \to \infty$.
\end{thm}

\proof  Using  \re{3.41} and  \re{3.6} together with
\begin{equation*}
\|v\|_{L^\infty_x}\leq\|v\|^{1/2}_{L^2_x} \|\partial_x v\|_{L^2_x}^{1/2}\label{3.43}
\end{equation*}
yields the desired result. 
\endproof

\section{Numerical Result}\label{sec6}
Numerical simulations are performed on several systems and  results on 
BBM-BBM and Bona-Smith systems  with
complete or partial dissipations are reported here. The results show  not only  that the 
theoretical results on the decay rates are sharp, but also the constants
involved are reasonably sized. 

In these numerical computations, the initial data are taken to be
\begin{equation*}
\begin{split}
\eta_0&=sech^2(\frac{\sqrt{2}}2(x-x_0)), \\
u_0&=\eta_0-\eta_0^2/4,
\end{split}
\end{equation*}
where $x_0$ is in the spatial domain $[0, L]$, where  $L$ is 
taken to be large enough so the solution near the boundary is smaller
than the machine roundoff error during the whole computation. 
The spectral method is used on the spatial domain $[0,L]$ and the leap-frog 
algorithm is used on the time advancing. 
The decay rate $r$ and the constant $C$  in 
$$\| v\| \sim C t^{-r}, \text{ as }  t \rightarrow \infty $$
is calculated by first
computing 
$$r(t_n):=-\frac{\log \frac{\| v \| (t_n)}{\| v\|(t_{n-1})}}{\log \frac{t_n}{t_{n-1}}}.
$$
The computation is stopped when $r(t_n)$ is approaching to a constant
and 
the value $r$ is  obtained 
by averaging the last 5 data. The constant $C$ is then
computed by averaging the last five $\|v\| (t_n) t_n^r$.

In the computations reported below,  $L=320$, $dx=0.1$ and
$dt=0.05$, where $dx$ and $dt$ are the meshsize in space and time 
respectively. 

\vskip .2cm 
\noindent
{\bf BBM-BBM system ($a=c=0, b=d=1/6$) with complete dissipation.} It is shown in Theorem
\ref{thm4.2} and \ref{thm5.1} that for {\it small data}, 
$$\| v\|_{L^2} \leq C_1 t^{-1/4} \quad \text{and} \quad \| v\|_{L^\infty} \leq C_2
t^{-1/2}.$$
The  numerical computation is performed for time
interval $[0, 50]$, and the result shows 
$$\| v\|_{L^2} \sim 1.4232 t^{-0.2470} \quad \text{and} \quad 
\| v\|_{L^\infty} \sim 1.4989
t^{-0.4963}.$$
Therefore, it is clear that   the theoretical result is sharp and the constants
involved are not large. Moreover, it seems that the small data
requirement might be removed if,  for example, other methods were
employed. 

\vskip .2cm 
\noindent
{\bf Bona-Smith system ($a=0, b=-c=d=1/3$) with complete and partial dissipation.}
This case is again covered by Theorem \ref{thm4.2} and \ref{thm5.1}. 
By direct computation, we obtain  for complete dissipation,  
$$\| v\|_{L^2} \sim 1.4015 t^{-0.2477}, \quad \| v\|_{L^\infty} \sim 1.4466
t^{-0.4998},$$
and for partial dissipation 
$$\| v\|_{L^2} \sim 0.6676  t^{-0.2519}, \quad \| v\|_{L^\infty} \sim 0.6595
t^{-0.5105}.$$

\vskip .2cm 
\centerline{Acknowledgment} Part of this work was performed when
the second author was enjoying the hospitality of the mathematics 
department in Purdue University.


\end{document}